\documentclass[letterpaper, 10pt, conference]{ieeeconf}  

\usepackage{amsmath}
\usepackage{amsthm}
\usepackage{amsfonts}
\usepackage{amssymb}
\usepackage{mathrsfs}
\usepackage{verbatim}
\usepackage{algorithm}
\usepackage{graphicx}
\usepackage{subfigure}
\usepackage{tabularx, booktabs}
\newcolumntype{Y}{>{\centering\arraybackslash}X}
\usepackage{multirow}
\usepackage{xcolor}

\newcommand{\R}{\mathbb{R}}

\newcommand{\calD}{\mathcal{D}}

\newcommand{\calY}{\mathcal{Y}}

\newcommand{\calT}{\mathcal{T}}

\newcommand{\calL}{\mathcal{L}}

\DeclareMathOperator*{\st}{\mbox{s.t.}}

\title{\LARGE \bf
Model Predictive Control of Discrete-Continuous Energy Systems \\ via Generalized Disjunctive Programming
}

\author{Arnab Bhattacharya$^{1}$, Xu Ma$^{1}$, Draguna Vrabie$^{1}$
}

\begin{document}
\maketitle
\thispagestyle{empty}
\pagestyle{empty}
\begin{abstract}
Generalized Disjunctive Programming (GDP) provides an alternative framework to model optimization problems with both discrete and continuous variables. The key idea behind GDP involves the use of logical disjunctions to represent discrete decisions in the continuous space, and logical propositions to denote algebraic constraints in the discrete space. Compared to traditional mixed-integer programming (MIP), the inherent logic structure in GDP yields tighter relaxations that are exploited by global branch and bound algorithms to improve solution quality. In this paper, we present a general GDP model for optimal control of hybrid systems that exhibit both discrete and continuous dynamics. Specifically, we use GDP to formulate a model predictive control (MPC) model for piecewise-affine systems with implicit switching logic. As an example, the GDP-based MPC approach is used as a supervisory control to improve energy efficiency in residential buildings with binary on/off, relay-based thermostats. A simulation study is used to demonstrate the validity of the proposed approach, and the improved solution quality compared to existing MIP-based control approaches.
\end{abstract}
\section{INTRODUCTION}
The past two decades have witnessed enormous interest in hybrid systems that combine continuous, time-driven dynamics with discrete, event-driven state transitions \cite{bemporad1999, borrelli2017}. The continuous dynamics are usually governed by a system of differential-algebraic equations or partial differential equations, while the discrete transitions are triggered by either state-dependent or time-dependent events \cite{barton94}. Such discrete-continuous behavior is observed in many existing energy systems. For example, electricity generators and demand-side loads drive the grid's continuous dynamics, while operation of protection devices or enforcement of controller hard limits exhibit discrete behavior \cite{esselman1999}. Similarly, HVAC scheduling in buildings involves switching HVAC equipments (e.g. boilers, chillers, fans etc.) on or off while optimizing  indoor-air temperatures with continuous physics-based dynamics \cite{risbeck2019}. 

While several studies exist on the analysis, verification and stability of hybrid systems \cite{doyen2018, minh13}, there remains significant challenges to control design, analysis and computation for hybrid systems \cite{blondel1999}. Moreover, there are very few approaches that can guarantee control performance in a computationally efficient way for general hybrid systems \cite{buss2002, attia2005}. An exception is mixed-integer linear programming (MILP), which has become increasingly popular for controlling hybrid systems with logic constructs, finite state-machines and linear discrete-time dynamics \cite{morari2001}. The recent success of MILP owes much to its flexibility, rigorousness, and extensive modeling capability, which is backed by significant improvements in computation speed of the state-of-the-art MILP solvers such as Gurobi and CPLEX. However, MILP models are NP-Complete, scale poorly with problem size, have high computational overhead and are unsuitable for systems with highly nonlinear dynamics. Moreover, MILP solution quality depends greatly on the strength of the model formulation used \cite{vielma2015}, which relies on a modeler's expertise to develop tractable and efficient formulations. Therefore, deployment of MILP-based control strategies for realistic hybrid systems remains a major challenge. 

Generalized disjunctive programming (GDP) provides an alternative to MILP for modeling systems with discrete and continuous components \cite{RAMANgrossman94, castro2012}. Specifically, a GDP model includes Boolean and continuous variables that are specified using mixed-integer constraints, logic disjunctions and logic propositions \cite{grossman2012pub}. Compared to traditional MILP, the logic structure in GDP can be exploited to reduce combinatorics, yield tighter relaxations, and improve bounds in branch and bound algorithms \cite{ruiz2017}. From a practical perspective, GDP facilitates modular integration of logic and algebraic constraints in an unified way, requires no modeling expertise, and produces models that are easy to interpret. 

In this paper, we provide a general framework for optimal control of hybrid systems using GDP. Specifically, we propose a GDP-based model predictive control (MPC) model for  piecewise-affine (PWA) systems with implicit switching logic. The proposed MPC model can be easily extended to other hybrid systems with complex logic rules and dynamical behavior. Compared to other MILP-based approaches, the proposed MPC model is easy to interpret and yields tighter relaxations that improve solution quality. We use the proposed MPC model to improve economics of conventional relay-based thermostats in residential buildings. Finally, we validate the performance and evaluate the improvement in solution quality compared to other MILP-based MPC approaches. 

The remainder of the paper is organized as follows. The next section introduces the general structure of a GDP model and related model reformulations. In Section \ref{sec:gdpmpc}, we present a general GDP-based MPC model for a general class of hybrid systems. Section \ref{sec:thermo_problem} presents the thermostat control problem for residential buildings. Section \ref{sec:numerical} illustrates the main results of our simulation study. Finally, we provide a few concluding remarks and directions for future work.
\section{Generalized Disjunctive Programming}\label{sec:gdp}
\noindent Generalized Disjunctive Programming (GDP) provides a high-level framework to model optimization problems with conditional logic constraints. A typical GDP model comprises of discrete and continuous variables that are related via algebraic constraints, logical disjunctions, and logic propositions \cite{ruiz2017} as follows:
\begin{subequations}\label{gdp}
\begin{align}
    z_\pi = \,\,\min_{y\in\calY} &\quad V^\pi(y) + \sum_{i=1}^L\gamma_i \label{gdp_obj} \\ 
    \mbox{s.t.}     &\quad g^\pi(y) \le 0,\label{gdp_global}\\
     &\quad\bigvee_{i=1}^L\begin{bmatrix}
                    a_i \\ 
                    r_i^\pi(y)\le 0\\
                    \gamma_i = c_i
                    \end{bmatrix} \label{gdp_disjunction} = \mathrm{True},\\ 
                    &\quad \Omega(a_1,\ldots, a_L)=\mathrm{True}, \label{gdp_intlogic}\\
                    &\quad a_1, \ldots, a_L \in \{\mathrm{True, False}\},
\end{align}
\end{subequations}
where $\pi$ is a parameter vector, $z_\pi\in\R$ is the optimal cost for a given $\pi$, $y\in \calY\subseteq \R^p$ is the decision vector, $g^\pi(y) \le 0$ is a set of global constraints, and $\{a_i: i=1,\ldots, L\}$ is a set of logic variables that can be either true or false. When $a_i$ is true, the local constraints in the set 
\[\calD_i\equiv\begin{bmatrix}
                    a_i \\ 
                    r_i^\pi(y)\le 0\\
                    \gamma_i = c_i
                    \end{bmatrix},\,\,\,i\in\{1, \ldots,L\},\]
are imposed, else they are relaxed. We refer to $\calD_i$ as the $i$-th disjunctive set associated with the logic variable $a_i$. The disjunction (or the logical OR) operator $\bigvee$ in \eqref{gdp_disjunction} ensures that only one of the disjunctive sets is active during optimization. 
We assume that $g^\pi(y)$ and $r_i^\pi(y)$ are convex in $y$ in the rest of the paper; see \cite{ruiz2017} for the nonconvex case. The variable $\gamma_i\in\R$ denotes the fixed cost of selecting the $i$-th disjunctive set. The logic proposition, $\Omega(a_1,\ldots, a_L)$, in \eqref{gdp_intlogic} models the relationships between the logic variables expressed in the Conjunctive Normal Form \cite{chandru1999}. The objective function in \eqref{gdp_obj} comprises of a variable cost $J^\pi(y)\in\R$ and a fixed cost of selecting a disjunctive set. Note that the feasible set in \eqref{gdp} is nonconvex as the union of non-overlapping disjunctive sets is nonconvex even if they are individually convex. To solve \eqref{gdp} as a mathematical program, the disjunctive and the logic constraints in \eqref{gdp_disjunction} and \eqref{gdp_intlogic}, respectively, are recast into equivalent algebraic constraints using binary (0-1) variables. Next, we discuss two well-known equivalent formulations of \eqref{gdp} that is amenable to solution via MIP solvers.


\subsection{Big-M Formulation}\label{sec:bigM}
\noindent The big-M method \cite{nemhauser1988} transforms each local constraint in a disjunctive set into an equivalent mixed-integer constraint 
\begin{equation}\label{bigMcons}
r_i^\pi(y)\le M_i(1-s_i),\quad i\in\{1,\ldots,L\},
\end{equation}
where $s_i\in\{0,1\}$ is defined such that $s_i=1$ when $a_i$ is true, and $s_i=0$ otherwise. The parameter $M_i\in(0,\infty)$ satisfies the condition $M_i\geq \sup\{r_i^\pi(y):\}_{y\in\calY}$ and is usually set to a large value for practical purposes. Define $\calL\equiv\{1,\ldots, L\}$. Then, it is straightforward to show that \eqref{gdp} is equivalent to the following MIP model:

\begin{subequations}\label{gdp_bigM}
\begin{align}
    z_\pi = \,\,\min &\quad V^\pi(y) + \sum_{i\in\calL}c_is_i \\ 
    \mbox{s.t.}     &\quad g^\pi(y) \le 0,\\
                    &\quad r_i^\pi(y)\le M_i(1-s_i),\quad \forall i\in\calL,\\ 
                    &\quad\sum_{i\in\calL}s_i=1,\\
                    &\quad As \geq a, \label{gdp_bigM_intlogic}\\
                    &\quad y\in\calY, s_i\in\{0, 1\}, \quad \forall i\in\calL.
\end{align}
\end{subequations}
Note that the constraint in \eqref{gdp_bigM_intlogic} is the algebraic representation of the logic proposition in \eqref{gdp_intlogic}. We refer to \eqref{gdp_bigM} as the \textit{big-M formulation}.
\subsection{Convex-Hull Formulation}\label{sec:chull}
\noindent The convex-hull formulation is based on the key result that the convex hull of the union of the sets $\calD_i$ can be expressed using the perspective functions of  $r_i^\pi(\cdot)$ and a set of disaggregated variables $\{y_i\}_{i\in\calL}$ defined for each $i\in\calL$ \cite{balas1985, raman1994}. The convex-hull formulation of \eqref{gdp} is an MIP model of the form (see \cite{raman1994} for more details):
\begin{subequations}\label{gdp_chull}
\begin{align}
    z_\pi = \,\,\min &\quad V^\pi(y) + \sum_{i\in\calL}c_is_i \\ 
    \mbox{s.t.}     &\quad g^\pi(y) \le 0,\\
                    &\quad y = \sum_{i\in\calL} y_i, \\
                    &\quad\sum_{i\in\calL}s_i=1,\\
                    &\quad s_ir_i^\pi(y_i/s_i)\le 0,\quad \forall i\in\calL,\label{perspective}\\ 
                    &\quad \ell s_i \leq y_i \leq us_i, \quad\,\,\, \forall i\in\calL,\\
                    &\quad As \geq a,\\
                    &\quad y\in\calY, s_i\in\{0, 1\}, \,\,\, \forall i\in\calL,
\end{align}
\end{subequations}
where $\ell$ and $u$ are lower and upper bounds of $y$. Note that the left-hand side of the constraint in \eqref{perspective} is the perspective function of $r_i^\pi(\cdot)$, which is convex when $r_i^\pi(\cdot)$ is convex \cite{soares1999}.



\begin{figure}[htb!]
    \centering
    \includegraphics[scale=0.58]{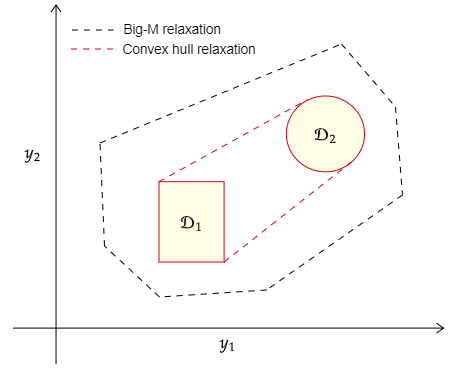}
    \caption{The convex hull model produces tighter continuous relaxations than the big-M model, resulting in improved branch and bound performance.}
    \label{fig:disjunction}
\end{figure}

The big-M and convex-hull formulations have complementary strengths and weaknesses. While the big-M formulation is straightforward and requires fewer variables, it usually produces weaker relaxations within branch-and-bound algorithm, especially for large values of $M_i$ \cite{lodi2010}. By comparison, the convex-hull model produces tighter relaxations that improves solution quality and convergence rate \cite{ruiz2017, grossman2003} but has higher computational overhead due to a larger number of variables and constraints in the model. Figure \ref{fig:disjunction} depicts the different relaxations obtained from the big-M and convex-hull formulations.
Next, we present a general GDP-based MPC model for hybrid systems.
\section{Disjunctive Model Predictive Control}\label{sec:gdpmpc}
\noindent In this section, we focus on piecewise-affine (PWA) systems as they are sufficiently expressive to model a large variety of hybrid processes with a high degree of accuracy \cite{bemporad2000, moehle2015perspective}. PWA systems are described by the state-space equations
\begin{subequations}\label{lti}
\begin{align}
    x_{t+1} &= A_{\lambda_t}x_t + B_{\lambda_t}u_t + E_{\lambda_t}d_t,\\ 
        y_t &= C_{\lambda_t}x_t + F_{\lambda_t}w_t,
\end{align}
\end{subequations}
where $x_t\in\R^{n}$ is the system state, $u_t\in\R^{m}$ is the control input, $y_t\in\R^k$ is the output,  $\lambda_t\in \Lambda\equiv\{1,\ldots, S\}$ is a switching regime, $S$ is the number of switching regimes, and $d_t\in\R^{v}$ and $w_t\in\R^{h}$ are noise terms at time $t\in\calT$, where $\calT$ is a finite planning horizon. The matrices $A_i\in\R^{n\times n}, B_i\in\R^{n\times m}$, $E_i\in\R^{n\times v}$, $C_i\in \R^{k\times n}$, and $F_i\in \R^{k\times h}$ are defined for each  $i\in\Lambda$. For each $t\in\calT$ and $i\in\Lambda$, define a binary variable $w_t^i\in\{0,1\}$ such that
\begin{equation*}
    w_t^i = 
    \begin{cases}
    1, &\textrm{if } \lambda_t = i,\\
    0, &\textrm{otherwise.}
    \end{cases}
\end{equation*}
Whenever the system is in regime $i$, the following set of constraints hold:
\[(w_t^i = 1) \implies \begin{cases}
x_{t+1}= A_ix_t + B_iu_t + E_id_t,\\
y_t= C_ix_t + F_iw_t,\\
r_i(x_t, u_t) \le 0,
\end{cases}\]
where $r_i(x_t, u_t)\le 0$ is a set of constraints local to regime $i$. Define $w_t \equiv (w_t^i)_{i\in\Lambda}$. The discrete dynamics are modeled using an implicit switching function
\begin{equation}\label{switchingfunc}
    w_{t+1} = \phi(w_t, x_t, u_t), \quad t\in\calT,
\end{equation}
where future regimes depend on the current regime, state and input. Assume that $\lambda_0$ (and therefore, $w_0$) is known with certainty. Note that $\phi$
may comprise of algebraic equations, logic propositions or both, which can be transformed into mixed-integer constraints. With a slight abuse of notation, define disjunctive sets for each $i\in\Lambda$ and $t\in\calT$ as follows: 
\[\calD_t^i = \bigvee_{i=1}^L\begin{bmatrix}
                    w_t^i \\ 
                    x_{t+1}= A_ix_t + B_i u_t + E_i d_t\\
                    y_t= C_ix_t + F_i w_t\\
                    r_i(x_t, u_t) \leq 0 \\
                    \gamma_t = f_i(y_t, u_t)
                    \end{bmatrix},\]
where $f_i(y_t, u_t)$ is assumed to be convex and represents the cost incurred at time $t$ in regime $i$. The disjunctive MPC model for the PWA system in \eqref{lti} is then expressed as

\begin{subequations}\label{gdp_mpc}
\begin{align}
\min & \quad \sum_{t\in \calT} \gamma_t \\
\st  & \quad x_{\min}\le x_t \le x_{\max}\qquad\qquad\qquad\quad\,\,\,\forall t\in\calT, \\
     &\quad u_{\min}\le u_t \le u_{\max} \qquad\qquad\qquad\quad\,\,\,\forall t\in\calT,\\
     &\quad \bigvee_{i=1}^L\calD_t^i = \mathrm{True}\qquad\qquad\qquad\qquad\quad\forall t\in\calT,\\
    &\quad \sum_{i=1}^L w_t^i = 1\qquad\qquad\qquad\qquad\qquad\,\,\,\,\forall t\in\calT,\\
    &\quad w_{t+1} = \phi(w_t, x_t, u_t) \qquad\qquad\qquad\,\,\,\,\forall t\in\calT,\\
    &\quad x_t\in\R^n, u_t\in \R^m, w_t \in \{0, 1\}^S \qquad\forall t\in\calT.
\end{align}
\end{subequations}
In the next section, we show how the disjunctive MPC model can be applied to a real-life problem in energy systems.   
\section{Optimal Control of Thermostats in Buildings}\label{sec:thermo_problem}
\subsection{Preliminaries and Modeling Assumptions}
\begin{figure}[htb!]
    \centering
    \includegraphics[scale=0.75]{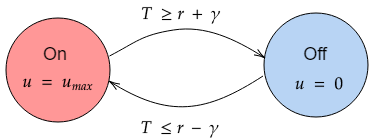}
    \caption{Transitions between the thermostat's operating states occur whenever the indoor temperature $T$ lie outside the switching bounds, i.e., $T\not\in[r-\gamma, r+\gamma]$, which leads to either maximum heating ($u=u_{\max}$) when the thermostat is On or no heating ($u=0$) when it is Off.}
    \label{fig:thermo_mode}
\end{figure}

\noindent We consider the problem of minimizing energy consumption and thermal discomfort in a single-zone residential building whose indoor temperature is controlled by a relay-based thermostat, as described in \cite{drgovna2015mpc}. The thermostat has two operating states, called the On and Off states. The thermostat prescribes maximum heating in the On state, while there is no heating in the Off state. Switching is triggered whenever the indoor temperature ($T_t$) violates thermostat-specific upper and lower switching bounds (known \textit{a priori}) around the reference setpoint ($r_t$), as illustrated in Figure \ref{fig:thermo_mode}. Specifically, the thermostat switches from the On to the Off state if $T_t > r_t + \gamma$, while the opposite occurs when $T_t < r_t-\gamma$, where $\gamma>0$ denotes half the width of the switching interval. MPC is used to optimize the reference setpoint over a finite planning horizon while enforcing the thermal comfort bounds, zone dynamical equations, and thermostat switching logic. We refer the reader to \cite{drgovna2015mpc} for more details on the problem setup.     

\noindent Before introducing the MPC model, we state the following key assumptions:
\begin{enumerate}
    \item The thermal dynamics in the zone are modeled as a linear time-invariant system of the form
    \begin{subequations}\label{lti_building}
    \begin{align}
        x_{t+1}&= Ax_t + Bu_t + Ed_t,\\
        T_t    &=Cx_t,
    \end{align}
    \end{subequations}
    where $x_t\in\R^n$ is the system state, $u_t\in\{0, u_{\max}\}$ is the (binary) heating control input, $d_t\in\R^v$ is the vector of exogenous disturbances (ambient temperature, indoor occupancy etc.), and $T_t$ is the indoor-air temperature at time $t\in\calT$.
    \item The system matrices $A, B, C$ and $E$ in \eqref{lti_building} are independent of the thermostat's operating states.
    \item The vectors $x_t$ and $d_t$ can be measured or estimated.
    \item Only thermal comfort is considered. Specifically, the indoor temperature $T_t$ must satisfy the constraint 
    \begin{equation}\label{comfort_orig}
        w - \theta \le T_t \le w + \theta, \quad t\in\calT,
    \end{equation}
    where $w$ is a user-defined temperature setpoint and $\theta>0$ is a user-defined band for thermal comfort (which can be different from $\gamma$). The terms $w+\theta$ and $w-\theta$ are the upper and lower comfort bounds, respectively, which are kept fixed over the planning horizon.
\end{enumerate}
\subsection{Model Predictive Control Model}
\noindent Using the current measurements of the state and disturbances, an optimal sequence of reference setpoints is obtained by solving the following MPC model over a finite prediction horizon $\calT$:
\begin{subequations} \label{mpc_thermo}
\begin{align}
\min\quad & \sum_{t\in\calT} (\alpha u_{t} + \beta m_t) \label{mpc_obj}\\
\st \quad  & x_{t+1}=Ax_{t}+Bu_{t}+Ed_{t} \qquad \quad\,\,\,\forall t\in\calT, \label{mpc_dynamics}\\
& T_{t}=Cx_{t} \qquad\qquad\qquad\qquad\qquad\,\,\forall t\in\calT, \\
& w-\theta -m_t \leq T_t \leq w+\theta +m_t  \quad\forall t\in\calT,\label{mpc_comfort}\\
& u_t=\left\{\begin{array}{ll}{u_{\max }} & {\text {if } s_t=1} \\ 
                               {0}        & {\text {if } s_{t}=0}
                               \end{array}\right. \qquad\quad\,\,\,\forall t\in\calT, \label{mpc_indicator} \\
& s_{t+1}=\phi(s_t, T_t, r_t) \quad\forall t\in\calT, \label{mpc_switch_func} \\
& r_t, T_t\in \R, m_t\in\R_+, x_t\in\R^n\qquad\,\,\,\forall t\in\calT, \\
& s_t\in\{0,1\} \qquad\qquad\qquad\qquad\qquad\,\,\forall t\in\calT.
\end{align}
\end{subequations}
The objective function in \eqref{mpc_obj} minimizes the total heating energy and thermal discomfort, where $\alpha$ and $\beta$ are weight parameters, and $m_t$ is a slack variable that measures deviation from the comfort bounds in \eqref{comfort_orig}. The constraint in \eqref{mpc_indicator} maps the heating control input to the operating state of the thermostat using a binary variable $s_t$ defined as
\begin{equation}\label{indicator_var}
    s_t = \begin{cases}1, &\textrm{if thermostat is currently in the On state},\\
                      0, &\textrm{otherwise}.
\end{cases}
\end{equation}
Moreover, the constraint in \eqref{mpc_switch_func} uses the thermostat's switching logic to map the operating states between consecutive periods using a switching function of the form
\begin{align*}
    \phi(s, T, r)=\left\{\begin{aligned} 1, & \quad \text { if } (s=1 \ \wedge \sim(T\geq r+\gamma)) \\ 
    & \quad \quad \vee\ (s=0 \wedge (T\leq r-\gamma)),\\ 
                                      0, & \quad \text { otherwise. } \end{aligned}\right.
\end{align*}
Using the big-M method and propositional calculus \cite{chandru1999}, the authors in \cite{drgovna2015mpc} transformed the MPC model in \eqref{mpc_thermo} into an MILP model. Next, we present an alternative GDP-based formulation of the MPC model in \eqref{mpc_thermo} that is amenable to the solution methods described in Section \ref{sec:gdp}.

\subsection{Model Reformulation via Disjunctive Programming}
\noindent The key idea here involves defining a finite collection of mutually exclusive disjunctive sets, where each set contains a group of MPC constraints that characterize a distinct operating mode of the system. In our context, an operating mode is defined by aggregating a possible combination of the thermostat's operating states (On or Off) without violating the switching logic. Table \ref{op_mode} lists four possible operating modes - by aggregating operating states in successive periods, denoted by $s_t$ and $s_{t+1}$, respectively - and the group of constraints imposed in each mode. For example, Mode 1 is encountered when the thermostat is On during successive periods, i.e., $s_t=1$ and $s_{t+1}=1$, which implies that the indoor temperature must be below the upper switching bound, i.e., $T_t < r_t+\gamma$, and maximum heating occurs in both periods, i.e., $u_t = u_{t+1}=u_{\max}$; all other operating modes have similar interetations. 
\begin{table}[htb!]
\centering
\caption{Disjunctive Operating Modes of the System.}\label{op_mode}
\begin{tabular}{l|l|l|l|}
\cline{1-4}
\multicolumn{1}{|c|}{Mode} & $s_t$  & $s_{t+1}$ &Local constraints \\ \hline
\multicolumn{1}{|c|}{1} & 1 (On)    & 1 (On)            & $T_t < r_t + \gamma$, $u_t=u_{t+1}=u_{\max}$      \\ \hline
\multicolumn{1}{|c|}{2} & 1 (On)    & 0 (Off)           & $T_t \geq r_t + \gamma$, $u_t=u_{\max}, u_{t+1}=0$ \\ \hline
\multicolumn{1}{|c|}{3} & 0 (Off)   & 1 (On)            & $T_t \leq r_t - \gamma$, $u_t=0, u_{t+1}=u_{\max}$  \\ \hline
\multicolumn{1}{|c|}{4} & 0 (Off)   & 0 (Off)           & $T_t > r_t - \gamma$, $u_t=u_{t+1}=0$ \\ \hline
\end{tabular}
\end{table}
Define a binary variable
\begin{equation}
    w_t^i = \begin{cases}
    1, &\textrm{if Mode $i$ is encountered at time $t$},\\
    0,  &\textrm{otherwise}.
    \end{cases}
\end{equation}
for each mode $i\in\{1, 2, 3, 4\}$ at time $t$. Then, the corresponding disjunctive sets are of the form


\begin{align}
\calD_t^1 = 
\begin{bmatrix}
w_t^1 \\
u_t = u_{\max} \\
u_{t+1} = u_{\max} \\
T_t < r_t + \gamma
\end{bmatrix};
&\quad \calD_t^2 = 
\begin{bmatrix}
w_t^2 \\
u_t = u_{\max} \\
u_{t+1} = 0 \\
T_t \geq r_t + \gamma
\end{bmatrix};\\
\calD_t^3 = 
\begin{bmatrix}
w_t^3 \\
u_t = 0 \\
u_{t+1} = u_{\max}\\
T_t \leq r_t - \gamma
\end{bmatrix};
&\quad \calD_t^4 =
\begin{bmatrix}
w_t^4 \\
u_t = 0 \\
u_{t+1} = 0 \\
T_t > r_t - \gamma
\end{bmatrix}.
\end{align}
The GDP-based MPC formulation is then expressed as 
\begin{subequations}\label{mpc_gdp}
\begin{align}
\min\quad &  \sum_{t\in\calT} (\alpha u_{t} + \beta m_t)  \\
\st \quad & \textrm{constraints \eqref{mpc_dynamics} -- \eqref{mpc_comfort}} \qquad\qquad\,\,\forall t\in\calT,\\
& \sum_{i=1}^4 w_t^i = 1 \qquad\qquad\qquad\qquad\qquad\forall t\in\calT,\\
&\calD_t^1\vee\calD_t^2\vee\calD_t^3\vee\calD_t^4 = \mathrm{True}\qquad\quad\forall t\in\calT, \\
& r_t, T_t\in \R, m_t\in\R_+, x_t\in\R^n\qquad\,\,\,\forall t\in\calT, \\
& w_t^1, w_t^2, w_t^3, w_t^4\in\{0,1\} \qquad\qquad\quad\forall t\in\calT.
\end{align}
\end{subequations}

\noindent The GDP model in \eqref{mpc_gdp} offers two key advantages over the MIP model proposed in \cite{drgovna2015mpc}. First, from a modeling perspective, the GDP model is significantly more interpretable as it avoids cumbersome algebraic transformations for the logic constraints in \eqref{mpc_indicator} and \eqref{mpc_switch_func}. Specifically, the  binary variables in the GDP model indicate the occurrence (or not) of operating modes with distinct characteristics. By comparison, the algebraic transformations in the MIP model results in auxiliary variables that have no straightforward interpretation. Second, the GDP model provides higher solution quality due to the tighter relaxations of the convex-hull formulation introduced in Section \ref{sec:chull}. 
\section{Simulation Experiments}\label{sec:numerical}
\subsection{Experimental Setup}
\noindent Similar to \cite{drgovna2015mpc}, we considered four state variables: floor temperature, internal facade temperature, external facade temperature and indoor temperature. Using a sampling time of 0.25 minutes ($\approx$ 15 seconds), we used the data from \cite{van2005integrated}  to estimate the matrices $A$, $B$, $C$ and $E$ in \eqref{lti_building}:
  \[A =10^{-2}\cdot\begin{bmatrix}
    99.97 & 0.00    & 0.00    & 0.00 \\
    0.00  & 99.98   & 0.00    & 0.00 \\
    0.00  & 0.00    & 99.92   & 0.00 \\
    1.77  & 4.28    & 0.00    & 93.48
    \end{bmatrix};\,\,C = \begin{bmatrix} 0 \\ 0 \\ 0 \\ 1 \end{bmatrix}';\]
    \[B = 10^{-4}\cdot\begin{bmatrix}
    0.0001 \\ 0.0001 \\ 0.0000 \\ 0.4421
    \end{bmatrix};\,E = 10^{-2}\cdot\begin{bmatrix}
    0.00    & 0.00    & 0.00 \\
    0.00    & 0.00    & 0.00 \\
    0.08    & 0.00    & 0.00 \\
    0.47    & 0.00    & 0.00
    \end{bmatrix}.\]
For simplicity, we neglect the exogenous disturbances, i.e., $d_t=0$ in \eqref{lti_building}. The user-defined temperature setpoint $w$ is kept fixed at $21^\circ C$ and the thermal comfort band $\theta$ is set to $1^\circ C$. The thermostat is assumed to be in the Off state initially, i.e., $s_0=0$. The maximum heating power $u_{\max}$ is set to 4 kW, while the width of the thermostat's switching bounds $\gamma$ is fixed at $1^{\circ} C$. We set $\alpha = 1$ and $\beta=10^5$ in the MPC objective function. The initial system state is set to $x_0=[21 \ 21 \ 21 \ 21]^\top$ and the starting time of all simulations is 7 A.M. Both the disjunctive-MPC model in \eqref{mpc_gdp} and the traditional MPC model of \cite{drgovna2015mpc}  were coded in YALMIP and solved using the Gurobi MILP solver on a 64-bit, Intel Core i7, 5th Gen., 16GB, 2.9GHz Windows machine. 

\subsection{Disjunctive MPC versus Relay-Based Thermostat Control}
\noindent First, we compare the disjunctive MPC (D-MPC) strategy to the pure relay-based thermostat control (RTC) in terms of the total energy consumption over a simulated horizon. We use a 2 hour horizon with 480 sampling periods. For the RTC strategy, we set $r = w = 21^\circ C$ for the entire horizon. In what follows, let $N$ be the number of MPC prediction periods and $M$ be the number of periods between successive MPC evaluations; it is assumed that the thermostat setpoints in periods between two MPC evaluations are set to the last computed MPC setpoint. 

\begin{figure}[htb!]
    \centering
    \includegraphics[width=0.48\textwidth]{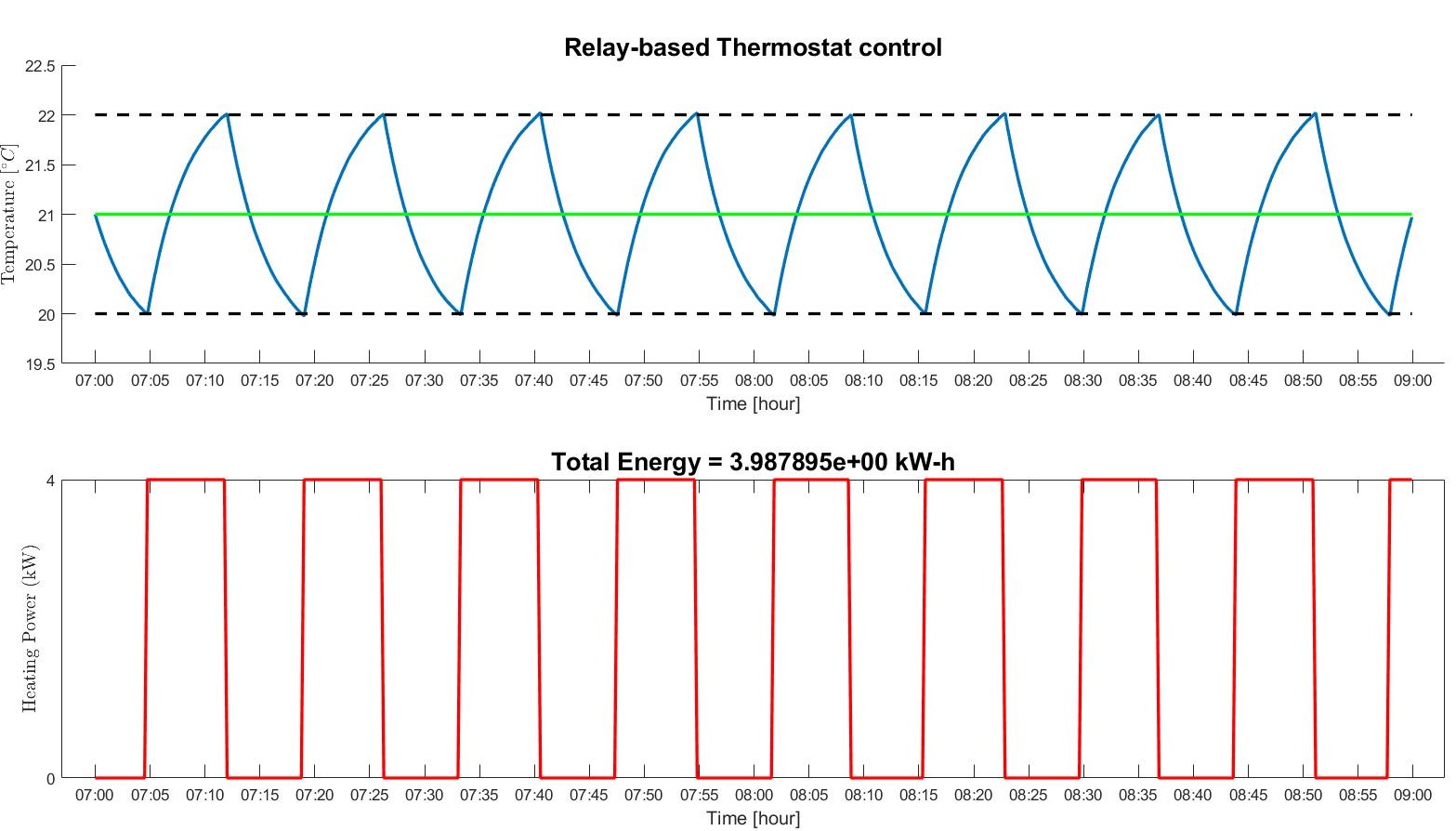}
    \caption{Temperature and heating profiles under relay-based thermostat control.}
    \label{fig:relay_control}
\end{figure}

\begin{figure}[htb!]
    \centering
    \includegraphics[width=0.48\textwidth]{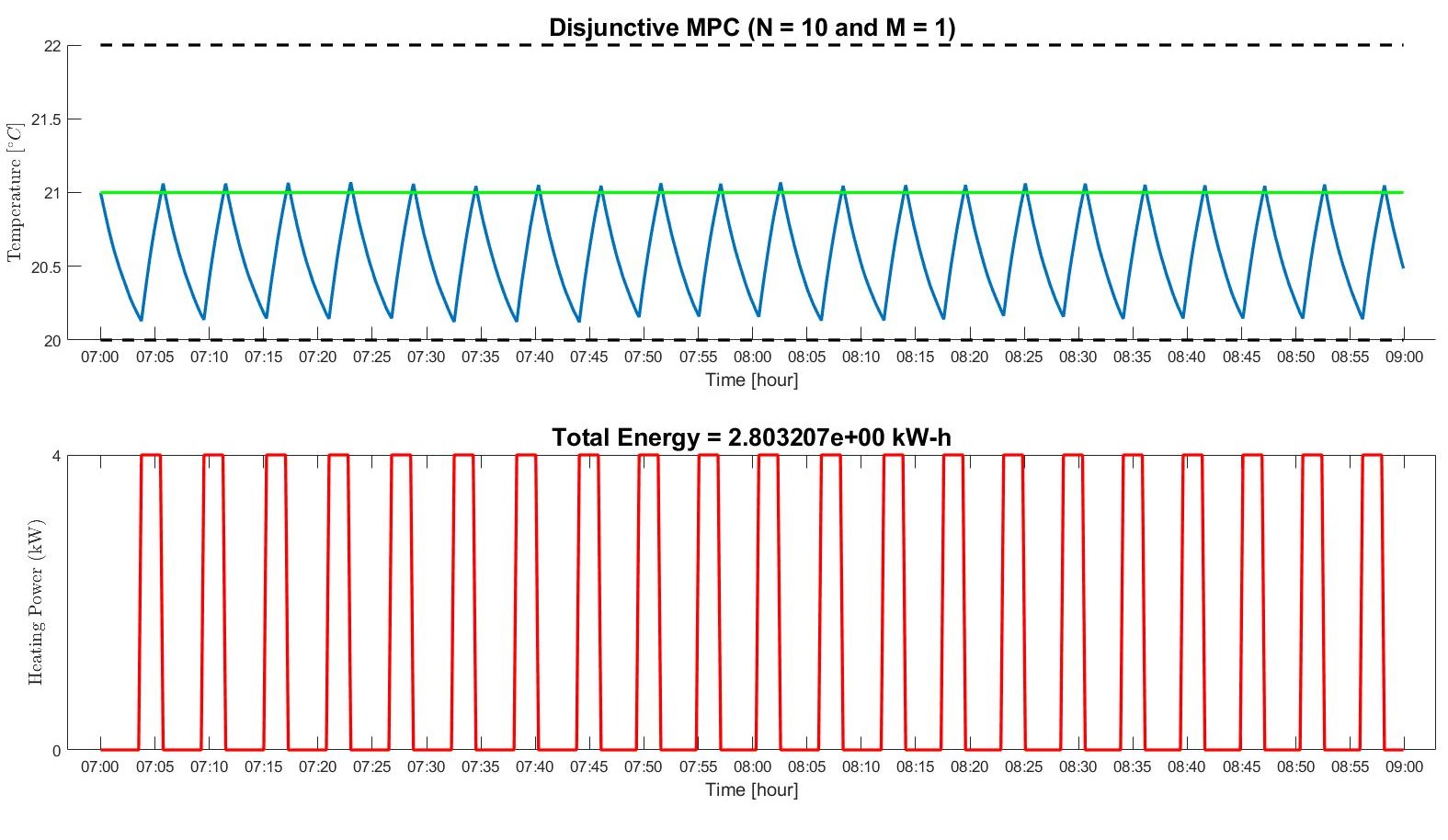}
    \caption{Temperature and heating profiles under disjunctive MPC with $N=10$ and $M=1$.}
    \label{fig:N10M1}
\end{figure}

Figure \ref{fig:relay_control} describes the zone temperature and heating profile under the RTC strategy. Both the plots are intuitive: the thermostat is switched on only when the temperatures go below the lower comfort bound ($20^{\circ} C$), and is switched off only when they exceed the upper comfort bound ($22^{\circ} C$). By comparison, Figure \ref{fig:N10M1} plots the temperature and heating profiles for the D-MPC strategy with $N=10$ and $M=1$ (MPC is solved at every sampling period). Note that D-MPC modulates the setpoints in a manner that switches the thermostat off more frequently and keeps the indoor temperature closer to the lower comfort bound for a longer duration of time. Consequently, D-MPC reduces energy consumption by around 30\% compared to RTC (from 3.98 kWh to 2.80 kWh). 

\begin{figure}[htb!]
    \centering
    \includegraphics[width=0.48\textwidth]{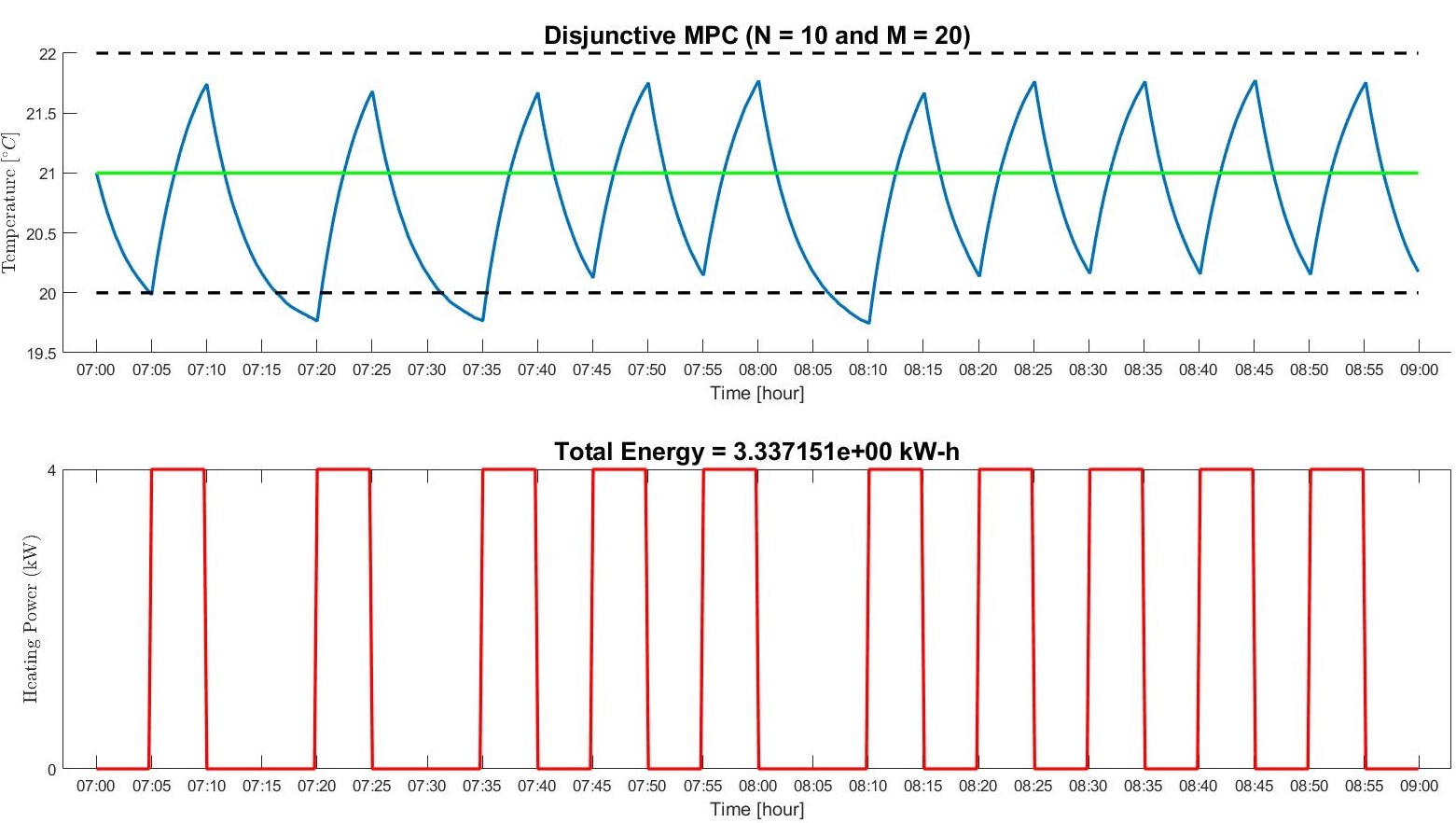}
    \caption{Temperature and heating profiles under disjunctive MPC with $N=10$ and $M=20$.}
    \label{fig:N10M20}
\end{figure}

\noindent However, choice of the MPC parameters can significantly impact control performance. For example, Figure \ref{fig:N10M20} plots the temperature and heating profiles for D-MPC when $N$ is still fixed at 10 but $M=20$ (MPC is solved every 20th sampling period). Reducing the frequency of MPC evaluations has a detrimental effect as energy consumption increased by 19\% compared to the case when $M=1$ (from 2.83 kWh to 3.37 kWh). This is because MPC's ability to react and take frequent corrective recourse reduces as $M$ increases. However, note that even for $M=20$, D-MPC reduced energy consumption by around 18\% compared to RTC (from 3.99 kWh to 3.34 kWh), which highlights the benefit of combining dynamic switching decisions with MPC's predictive capabilities.
\subsection{Disjunctive MPC versus Traditional MPC}
\begin{figure}[htb!]
    \centering
    \includegraphics[scale=0.47]{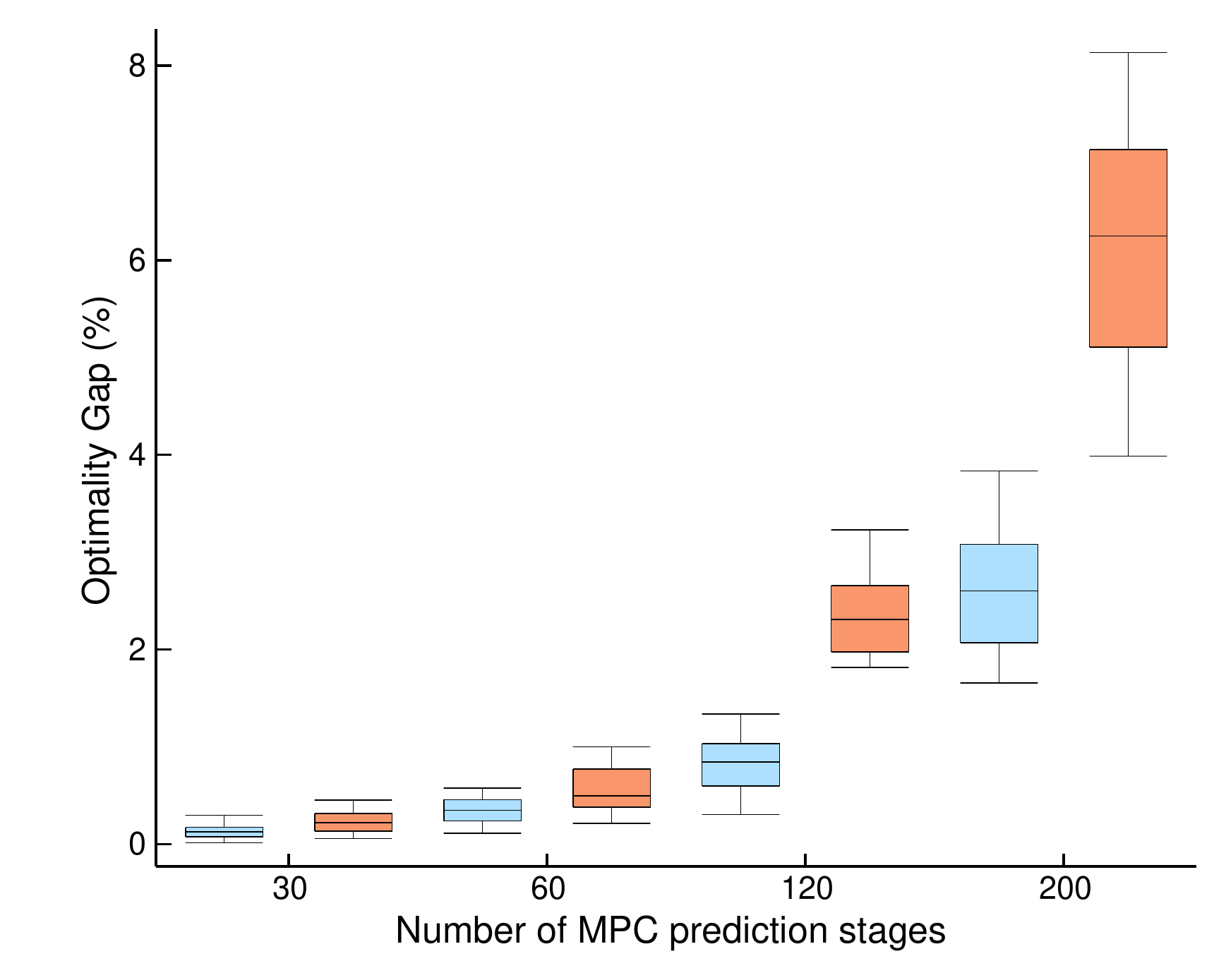}
    \caption{Comparison of the optimality gaps between D-MPC (blue) and MIL-MPC (orange) instances for different MPC prediction horizons.}
    \label{fig:optgap}
 \end{figure}

\noindent Next, we compare the solution quality of the D-MPC model with the MILP-based MPC model (MIL-MPC) of \cite{drgovna2015mpc}. We conducted an open-loop statistical study where 50 initial states $(x_0)$ were randomly sampled and solved for both MPC models for $N\in\{30,60,120,200\}$ ($50\times4=200$ instances for each MPC type). Because the problem sizes were not large, we solved all 400 instances to optimality (the big-M parameter was set to $10^4$ for MIL-MIP). Next, we restricted the number of iterations in the Gurobi solver to 30 and recorded the final solutions upon termination. For each MPC instance, we computed the optimality gap \% as
\[\textrm{Optimality gap \%} = \frac{z^*-\tilde{z}}{z^*}\times 100,\]
where $z^*$ is the optimal objective cost and $\tilde{z}$ is the objective cost at the end of the 30th iteration. Note that solution quality increases as optimality gap decreases. Figure \ref{fig:optgap} confirms that D-MPC produces solutions of higher quality (on average) compared to MIL-MPC due to the tighter continuous relaxations obtained from the convex hull formulation. Moreover, solution quality degrades rapidly for MIL-MPC for large problem instances. For example, when $N=200$, MIL-MIP has an average optimality gap of 5.91\% compared to 2.36\% for D-MPC. 

\section{CONCLUSIONS}
\noindent A GDP-based MPC model has been proposed for hybrid systems. Two well-known reformulation techniques for disjunctive constraints have been discussed, and their application to hybrid systems have been highlighted. A discrete-continuous energy system served as an illustration, and a simulation study validated the performance and improved solution quality of the proposed MPC model.





\bibliographystyle{IEEEtran}
\bibliography{root}

\end{document}